\input amstex
\input xy
\xyoption{all}
\input epsf
\documentstyle{amsppt}
\document
\magnification=1200
\NoBlackBoxes
\nologo
\hoffset=1.5cm
\voffset=1truein

\pageheight{16cm}

 \centerline{\bf FORGOTTEN MOTIVES:}
 
 \medskip
 
 \centerline{\bf THE VARIETIES OF SCIENTIFIC EXPERIENCE}
 \bigskip

\centerline{\ Yuri I. Manin}

\medskip

\centerline{\it Max--Planck--Institut f\"ur Mathematik, Bonn, Germany}
  
\bigskip

\hfill{\it Le gros public:\quad\quad\quad\quad\quad\quad\quad}

{\hfill A po\^ele, Descartes!  \`a po\^ele!\quad\quad\quad\quad}

\medskip

{\hfill(\it R.~Queneau, Les Oeuvres compl\`etes }

\hfill{\it  de Sally Mara)}

\bigskip

When I arrived in Bures--sur--Yvette in May 1967, the famous seminar SGA 1966--67, dedicated to the Riemann--Roch theorem,
was already drawing to an  end.
Mlle Rolland, then L\'eon Motchane's secretary at the IH\'ES, found for me
a nice small apartment in Orsay. Each early morning, awoken  to the loud chorus of singing birds,
I walked to Bures, anticipating the new session of private tutoring on the then brand--new project of motives,
by the Grand Ma{\^\i}tre himself, Alexandre Grothendieck. Several pages, written by his hand then,
survive in my archive;  in particular, the one dedicated to the ``Standard Conjectures''. These
conjectures remain unproved after 
half a century of vain efforts.  Grothendieck himself saw them as the cornerstone  of the whole project.
In the letter to me dated March 20, 1969, he wrote:

\smallskip

{\it Je dois avouer \`a ma honte que je ne sais plus distinguer \`a premi\`ere vue ce qui est
d\'emontrable (voire plus ou moins trivial) sans les conjectures
standard, et ce qui ne l'est pas. C'est \'evidemment honteux qu'on n'ait pas d\'emontr\'ees
les conjectures standard!}

\smallskip

Still, during the decades that have passed since then the vast realm of motives kept rewarding the humility of many researchers  prepared to be happy
with what they could do using the tools they could elaborate.

\smallskip

Several times Grothendieck invited me  to his house at rue de Moulon. He allowed me to browse through
his bookshelves; I borrowed a few books  to read at home. When I last visited him
a day or two before my departure, I asked him to sign  a book or  paper for me.
To my amazement, he opened ``Les {\OE}uvres compl\`etes  de Sally Mara'' by Raymond Queneau
and scribbled on the first page: 
\smallskip
\centerline{\it Hommage affectueux  R.~Queneau} 

\newpage

\centerline{\bf Early history of motives}

\medskip

Having returned to Moscow in June 1967, after five or six weeks of intense training
with Grothendieck, I spent several months writing down his main definitions related to motives and studying
necessary background material in the literature.  I was very pleased when it turned out that
I could answer one of his questions and calculate the motive of a blow--up
without using standard conjectures. My paper [Ma68] containing this exercise was submitted  next summer and published
in Russian. It became the first ever publication on motives, and
Grothendieck recommended it to David Mumford (in his letter of April 14, 1969)
as {\it  ``a nice foundational paper''} on motives.

\smallskip

 Grothendieck  wrote a letter  in Russian to me about this paper (05/02/1969).
This seems to be the only  document showing that he had some Russian, probably, learned from
his father.

\smallskip

The first step in the definition of a category of (pure) motives is this.
We keep objects of a given algebraic--geometric category, say of  smooth projective varieties over a fixed field $Var_k$,  
but replace its morphisms
 by {\it correspondences}.
This passage implies that morphisms $X\to Y$ now form an {\it additive group,}
or even a $K$--module rather than simply a set, where $K$ is an appropriate coefficient ring.
Moreover, correspondences themselves are not just cycles on $X\times Y$ but
{\it classes} of such cycles modulo an ``adequate'' equivalence relation. The coarsest such relation 
is that of numerical equivalence, when two equidimensional  cycles are equivalent if  their intersection indices
with each cycle of complementary dimension coincide. The finest one is the rational  (Chow) equivalence,
when equivalent cycles are deformations over a base which is a chain of rational curves.
Direct product of varieties induces tensor product structure on the category.

\smallskip

The second step in the definition of the relevant category of pure motives consists in 
a formal construction of new objects (and relevant morphisms) that are ``pieces'' of varieties: kernels and images of
projectors, i.~e. correspondences $p:\,X\to X$ with $p^2=p$. This produces
a {\it pseudo--abelian,}  or {\it  Karoubian} completion of the category. In this new category,
the projective line $\bold{P}^1$ becomes the direct sum of the (motive of) a point
and the Lefschetz motive $\bold{L}$ (intuitively corresponding to the affine line).

\smallskip

The third, and last step of the construction, is one more formal enhancement of the
class of objects: they now include {\it all} integer tensor powers  $\bold{L}^{\otimes n}$, not just
non--negative ones, and tensor products of these with other motives.

\smallskip

Various strands of intuition are interwoven in this fundamental pattern discovered by
Grothendieck, and I will now try to make them (more) explicit.

\smallskip

The basic intuition that guided Grothendieck himself, was the image of the category of pure Chow
motives $Mot_k$ as the receptacle of the ``{\it universal cohomology theory''} $V_k\to Mot_k$:\,
$V\mapsto h(V)$.
The universal theory was needed in order to unite various cohomological constructions,
 such as Betti, de Rham--Hodge, and \'etale cohomology. 
 
 \smallskip
 What looked paradoxical
 in this image was the following observation about transcendental cycles on an algebraic variety $X$.
 One could get hold of these cycles for $k=\bold{C}$ by appealing to algebraic topology,
 or else to complicated constructions of homological algebra involving all finite covers of $X$.
 
 \smallskip
 But in the category of pure motives,  from the start one dealt {\it only} with algebraic cycles,
 represented by correspondences, and it was intuitively not at all clear how on earth they could 
 convey information about
 transcendental cycles.  Indeed, the main function of the ``Standard Conjectures''  was to serve as
 a convenient bridge from algebraic to transcendental. Everything that one could prove without them
 was indeed {\it ``plus ou moins trivial''} --  until people started treating correspondences themselves
 using sophisticated homological algebra (partly generated by the development of  \'etale cohomology
 and Grothendieck--Verdier's introduction of derived and triangulated
 categories).
 
\smallskip

However, the passage from the set of morphisms to the $K$--module of correspondences
involves one more intuitive idea, and it can be most succinctly invoked
by referring to {\it physics}, namely the great leap from the  classical mode 
of description of nature to the quantum one. This leap defined the science of the XXth century.
Its basic and universal step consists in the introduction of a {\it linear span} of
everything that in classical physics was only a set: points of a phase space,
field configurations over a domain of space--time etc. Such {\it quantum superpositions}
then form  linear spaces on which Hilbert--like scalar products are defined,
that in turn allow one to speak about probability amplitudes, quantum observations etc.

\smallskip

I have no evidence that Grothendieck himself thought then about quantum physics 
in relation to his algebraic geometry project. We do know that concerns about
weapons of mass destruction and collaborationist behaviour of scientists towards their
governments and military--industrial complexes inspired in him deep disturbance and 
aversion.  The most direct source of his inspiration might have been algebraic
topology which, after the 1940s, laid more stress on chains and cochains than
on simplices and the ways they are glued together.

\smallskip

However, in my personal development  as a mathematician in the 1970's--80's and later,
the study of quantum field theory played a great role, and feedback from theoretical
physics  -- that was ahead --  to algebraic geometry  became  a great source of 
inspiration for me. I was and remain  possessed by a  Cartesian dream, poetic rationalism,
whatever  history has yet to say about {\it Der Untergang des Abendlandes.}
\smallskip

Below I will sketch a map of two branches of the  development of Grothendieck ideas about motives that
approximately followed two intuitions invoked above: from homological algebra and from physics respectively.
The references at the end of this essay constitute the bare minimum of the relevant research, but the reader
will be able to find much additional bibliographical material in the survey collection [Mo91] and in [A04], [VoSuFr00], [Ta11].

\bigskip

\centerline{\bf Motives and homological algebra}

\medskip

The most common linear objects are modules over rings in algebra and sheaves of modules
in algebraic geometry. {\it Free modules/locally free sheaves}  are the closest to classical linear spaces.
\smallskip
General algebraic variety $X$, or a scheme, is a highly non--linear object.

\smallskip
In  classical algebraic geometry over, say, the complex numbers, the variety $X$ used to be identified
with the topological space $X(\bold{C})$ of its $\bold{C}$--points, and could be studied
by topological methods involving triangulations or cell decompositions.
In the geometry over, say, finite fields, this did not work, and when in 1949 Andr\'e Weil 
stated his famous suggestion that point counting over finite fields should be done using 
trace of the Frobenius endomorphism acting upon appropriately defined  (co)homology groups of  $X$,
it generated a flow of research. 

\smallskip

The first product of this research was the creation of the cohomology theory
of coherent sheaves of modules $\Cal{F}$ on varieties $X$ or more general schemes.
Now, in a constructive definition of $H^*(X,\Cal{F})$,
one could either stress combinatorics of covers of $X$ by open sets in the Zariski topology ($\check{\roman{C}}$ech
cohomology)
or, alternatively, ``projective/injective resolutions'' of $\Cal{F}$, that is special exact complexes of sheaves\quad   
$\dots  \to F_2\to F_1\to F_0:=\Cal{F}\to 0$ or similarly with arrows inverted. This passage from the dependence of $H^*(X,\Cal{F})$ on  the non--linear argument $X$ to the dependence on  the linear argument  $\Cal{F}$ was very characteristic for the early
algebraic geometry of 1950's and 1960's. ``Homological Algebra'' by H.~Cartan and S.~Eilenberg,
the famous FAC, ``Faisceaux Alg\'ebriques Coh\'erents'' by J.-P.~Serre, became the standard handbooks
for every aspiring young algebraic geometer.

\smallskip

David Mumford and I started our training as algebraic geometers at the same time, about 1956, he at Harvard,
I at Moscow University. David reminisces that his teacher Zariski ``was motivated by the need to make the work of the Italian school rigorous by using the new methods of commutative algebra''.
My teacher Shafarevich also suggested to us  to study glorious Italian algebraic geometry,
approaching it armed with modern insights and techniques developed by  Serre, Grothendieck and their school. 

\smallskip

I had no time nor use for a course in  ``Instant Italian'', so I tried to read two books
simultaneously, ``Le Superficie Algebriche'' by Federigo Enriques (Zanichelli, 1949)
and ``La Divina Commedia'', and each time that I opened Enriques (or for that matter, SGA),
I recited mournfully:  {\it  lasciate ogni speranza voi ch'entrate ... }

\smallskip

Nevertheless, it worked. When I brought  xeroxed papers by Gino Fano back from Bures in 1967,
Vassya Iskovskikh and I could read them without bothering much in which language
they have been written, and then produce the first examples of  birationally rigid varieties,
and unirational but not rational threefolds
using Fano methods.

\smallskip

Homological algebra proved more resistant, and here I learned most of what I understand now from
the next generation of eager young Moscow students, who by now have been mature researchers themselves
for a long time.

\smallskip

We first learned, of course, about the basic Grothendieck--Verdier presentation of homological algebra as
the theory of derived, and more generally,  triangulated categories. Passage from the Bourbaki language of structures to the
now domineering language of categories (and then polycategories) involved several radical changes of
intuition, and as  is now clear, led into the garden of forking paths. The passage from one crossroad to 
another one always involved a decision about what should be disregarded, and later it could happen -- and always did happen --
that one was bound to turn back again and  recollect some forgotten ideas.

\smallskip

The story of {\it derived categories} started with categories, whose {\it objects} were complexes (of abelian groups/sheaves/objects of an abelian category)
considered modulo homotopy. 

\smallskip
In the framework of Grothendieck--Verdier  triangulated categories, one
forgot about initial objects--complexes and  focused on an abstract additive category, endowed with
a translation functor and a class of diagrams, called distinguished triangles.  But the problem of non--functoriality of 
cones led back to the complexes of abelian groups, this time upgraded to the level of {\it morphisms} rather
than objects. 
\smallskip
This was, of course, a special case of {\it enriched categories}, which in the simplest
incarnation postulate  Bourbaki--structured morphism sets $\roman{Hom}\,(X,Y)$, but with an upgrading: this time
one clearly had to deal with the case of  {\it categorified} morphism sets. However, when one allows  morphisms to be objects of a category,
then morphisms of this second floor category might form a category as well \dots and we find ourselves
ascending the Tower of Babel  that could  cause despair 
even in Grothendieck himself.

\smallskip

For the limited purposes of this note, I will disregard subtleties and various versions
of the notion of triangulated/dg--categories, and will only sketch several basic discoveries
of the last decades relating such categories with motives.

\smallskip

Roughly speaking, starting with a category of varieties (or  schemes) $X$, one may consider either the replacement of each $X$ by
 a triangulated category
$D(X)$ of complexes of (quasi)--coherent sheaves on $X$, or else return to the initial
Grothendieck insight, but replace correspondences by {\it complexes of correspondences.}

\smallskip
The latter approach led to the Voevodsky's motives ([VoSuFr00]). I will focus on some achievements of the first one.

\smallskip

One of the first great surprises was Alexander Beilinson's discovery ([Be83]) that a derived category
of a projective space can be described as a triangulated category made out of modules over a
Grassmann algebra. In particular, a projective space  became ``affine'' in some kind
of non--commutative geometry! The development of Beilinson's technique led to a general machinery describing triangulated categories in terms of exceptional systems and extending the realm of candidates to the role
of non--commutative motives. 
\smallskip

D.~Orlov ([Or05]) proved a general theorem to the effect that if $X$, $Y$ are smooth projective $k$--varieties
and if there is a fully faithful functor $F:\, D^b(X)\to D^b(Y)$, then the Chow motive $h(X)$ 
is a direct summand of  $h(Y)$ ``up to translations and twists by Lefschetz/Tate motives''.

\smallskip

M.~Kontsevich formalised the properties of $dg$--categories, expressing properness and smoothness
in case of the derived categories of varieties, and defined the respective class of categories (modulo homotopy) 
as  ``spaces'' in non--commutative algebraic geometry.  He then defined the respective class of Chow motives
and has shown that there exists a natural fully faithful functor embedding Grothendieck's Chow
motives (modulo twists) into non--commutative motives. These ideas were further developed by Tabuada, Marcolli, Cisinski et al.,
cf. the recent survey [Ta11] and references therein.

\bigskip

\centerline{\bf Motives and physics}

\medskip

In the mid--1970's and later, algebraic geometry interacted with physics more intensely that ever before:
self--dual gauge fields (instantons), completely integrable systems (Korteweg--de Vries equations),
emergence of supergeometry (based upon formal rules of Fermi statistics), the Mumford form and the Polyakov measure 
on moduli spaces of curves (quantum strings) have been discussed  at joint seminars
and local and international conferences of physicists and mathematicians.

\smallskip

Motives did not yet appear in this picture. However, in 1991 something new and unexpected happened.

\smallskip

B.~Greene in his book  ``The Elegant Universe. Superstrings,
Hidden Dimensions and the Quest for the Ultimate Theory'' tells the following story:

\smallskip

{\it ``At a meeting of physicists and mathematicians in Berkeley in 1991,
Candelas announced the result reached by his group using string theory and mirror symmetry:
317 206 375. Ellingsrood and Str{\o}mme announced the result of their
very difficult mathematical computation: 2 682 549 425. For days,
mathematicians and physicists debated: Who was right? [\dots ]

\smallskip

About a month later, an e--mail message was widely circulated
among participants in the Berkeley
meeting with the subject heading: {\it Physics Wins!}   Ellingsrood and Str{\o}mme had found
an error in their computer code that, when corrected, confirmed Candelas's result.''}

\smallskip

The problem about which Greene speaks is this. Consider a smooth hypersurface $V$ of degree 5
in $\bold{P}^4.$ Denote by $n(d)$ the (appropriately defined) number of rational curves of degree $d$
on $V$. Calculating $n(d)$ looks like perfectly classical problem of enumerative
algebra geometry, and in fact the numbers $n(1)=2875$ and $n(2)=609250$ were long known.
The physicists Ph.~Candelas,  X.~C.~de la Ossa, P.~S.~Green, and  L.~Parkes using machinery and heuristics of quantum string theory,
calculated not just $n(3)$, but gave an analytic expression for a total generating function
for these numbers, using the so called Mirror Conjecture. The mathematicians
G.~Ellingsrood and S.~A.~Str{\o}mme produced a computer code calculating $n(3)$.

\smallskip

Omitting a lot of exciting developments of this rich story, I will briefly explain only the part
that refers to the new and highly universal motivic structure that emerged in algebraic geometry.
I will speak about varieties, although in fact Deligne--Mumford stacks form the minimal habitat
for this structure, and the respective extension of the construction of pure motives
for them is needed; this was done by B.~To\"en.

\smallskip

Roughly speaking, we now treat the general problem, inherited from classical enumerative geometry: given a projective variety $V$, (define and) calculate the number of algebraic curves of genus $g$ on $V$,
satisfying additional incidence conditions that make this number finite, as in the Euclidean archetype:
"one line passes through two different points of plane''. After considerable efforts, one can
define for all stable values of $(g,n)$  a Chow class $I_{g,n}$ on $V^n\times \overline{M}_{g,n}$
with coefficients in the completed semi--group ring, say $\bold{Q}[[q^{\beta}]]$ where $\beta$ runs over
integral classes in the Mori cone of $V$. This class expresses the virtual incidence relation,
described above, by reducing it to the positions of the respective points in $V^n$
on the one hand, and to the position of the respective curve in the Deligne--Mumford stack
of curves of genus $g$ with $n$ marked points. 

\smallskip

When this is done, a list of universal properties of the classes  $I_{g,n}$ treated as motivic
morphisms, defines essentially the (co)action of the modular (co)operad with components
$h(\overline{M}_{g,n})$ in the category of motives upon each {\it total} motive $h(V)$
(I use the word total in order to stress that we are not allowed to pass to pieces here, although twisting and translations
are in fact present, cf. [BehM96]).

\smallskip

The sophistication of both theoretical (and imaginative) physics and abstract mathematics that
cooperated to discover this picture is really amazing, and I would like to draw attention to the fact that
our traditional (mis)representation of mathematics as a language and  technical tool needed to make
physical intuition precise, was reversed here: physical intuition helped discover
mathematical structures that were not known before.  One remarkable result of this was Deligne's generalisation
of the Tannakian Galois formalism ([De02]): it turned out that motivic Galois groups are actually supergroups,
so that the Fermi statistics now firmly resides in algebraic geometry as well, which
up to then was ``purely bosonian''.

\smallskip
Of course, such reversals have happened
many times in history, but here the contemporary status of both theory of motives and quantum strings
adds a strong romantic touch to the story. The beautiful  two--volume  cooperative project 
of the two communities trying to enlighten each other, [QFS99], is branded by two epigraphs.
The epigraph to the first volume is a quotation from Grothendieck's  ``R\'ecoltes et Semailles'':
\smallskip
{\it Passer de la m\'ecanique quantique de Newton \`a celle d'Einstein doit \^etre un peu,
pour le math\'ematicien, comme de passer du bon vieux dialecte proven{\c{c}}al \`a l'argot parisien
dernier cri. Par contre, passer \`a la m\'ecanique quantique, j'imagine, c'est passer du fran\c{c}ais
au chinois.}
\smallskip

(In the pre--post--modern times one would have said: ``It's all Greek to me!'').

\smallskip
The second volume starts with epigraph, written in Chinese logograms, from Confucius'  ``Analects'', 17:2.
Here I give its translation:

\smallskip

{\it The Master said: ``Men are close to one another by nature. They drift apart through behavior that is constantly repeated''.}
\smallskip

This is the  collective riposte of the two communities, 
 arguing their closeness, but in the language that is foreign for both.

\bigskip
\centerline{***}

\medskip

In his letter to me from Les Aumettes  dated  March 8, 1988, the last letter that I have,
Grothendieck has written:

\smallskip

{\it $ \dots$ thanks for your letter of birthday congratulations, and please excuse my being late in replying
to this letter, as well as the previous one and thanking for the reprint with dedication of november
last year. Your letter struck me as somewhat formal and kind of ill at ease, and surely
my silence has contributed to it. What I had to say about the
spirits in mathematics today I said in the volumes I sent you and a number
of other former friends. I am confident that before the year 2000 is reached, mathematicians
(and even non--mathematicians) will read it with care and be amazed about times strange
at last left behind $ \dots$}

\medskip

I met Grothendieck almost half a century ago.  Thinking back on his imprint on me then,
I realise that it was his generosity and his uncanny sense of humour that struck me most,
the carnivalistic streak in his nature, which I later learned to discern in other anarchists
and revolutionaries.

\medskip

On the front cover of  the issue no 14 of  {\it ``Survivre \dots et Vivre''  (Octobre--novembre 1972)}
that miraculously reached me by post in Moscow, I read:

\medskip

{\it 2 FRANCS

Canada 50 c

Communaut\'es:

1 fromage de ch\`evre.}

\vskip1cm

\centerline{\bf References}

\medskip

[Gr 69] A.~Grothendieck. {\it Standard conjectures on algebraic cycles.} In: Algebraic Geometry, Bombay
Colloquium, 1968, Oxford UP, 1969, pp. 193--199.

\smallskip

[Ma68] Yu. Manin. {\it Correspondences, motives and monoidal transformations.} Math. USSR--Sb., 6 (1968),
pp.439--470.

\smallskip

[Ja92] U.~Janssen. {\it Motives, numerical equivalence, and semi--simplicity.} Invent. Math. 107 (1992),
447--452.

\smallskip

[Mo91] {\it Motives.}  Proc. Summer Research Conference, 1991. Eds.  U.~Janssen, S.~Kleiman, J.--P.~Serre. Proc. Symp. in Pure Math., vols. 55.1, 55.2, AMS 1994.

\smallskip

[VoSuFr00] V.~Voevodsky, A.~Suslin, E.~Friedlander. {\it Cycles, transfers, and motivic homology theories.}
Ann.~Math.~studies, Princeton UP, 2000.
\smallskip
[A04] Y.~Andr\'e. {\it Une Introduction aux Motifs (motifs purs, motifs mixtes, p\'eriodes).} Panoramas
et Synth\'eses, Nr. 17, Soc. Math. de France, 2004.
\smallskip

[Or05] D.~Orlov. {\it Derived categories of coherent sheaves and motives.} arXiv:math/0512620

\smallskip

[Be83] A.~Beilinson. {\it The derived category of coherent sheaves on $\bold{P}^n$}. Selecta Math. Soviet. 3, 1983, pp. 233--237.
\smallskip
[Ta11] G.~Tabuada. {\it A guided Tour through the Garden of Noncommutative Motives.}
arXiv:1108.3787

\smallskip

[QFS99] {\it Quantum Fields and Strings: A Course for Mathematicians. Vols. 1, 2.} Ed. by P.~Deligne, Ed Witten, et al.,
AMS, 1999.

\smallskip

[De02] P.~Deligne. {\it Cat\'egories tensorielles.} (Dedicated to Yuri I. Manin on the occasion of his 65th birthday.) Mosc.~Math.~J. 2 (2002), no. 2, 227--248.
\smallskip

[BehM96] K.~Behrend, Yu.~Manin. {\it Stacks of  stable maps and Gromov--Witten invariants.}
Duke Math.~ J., 85, no.1, 1996, pp. 1--60.

\smallskip

\enddocument